\documentclass[10pt, a4paper]{article}
\usepackage{times}
\usepackage{hyperref}
\usepackage[british]{babel}
\usepackage{enumerate, longtable}
\usepackage{amsmath, amscd, amsfonts, amsthm, amssymb, latexsym, comment, stmaryrd, graphicx}
\usepackage[T1]{fontenc}
\usepackage[latin1]{inputenc}

\newcommand{\FF}{\mathbb{F}}
\newcommand{\NN}{\mathbb{N}}
\newcommand{\QQ}{\mathbb{Q}}
\newcommand{\ZZ}{\mathbb{Z}}

\newcommand{\Qbar}{\overline{\QQ}}
\newcommand{\Fbar}{\overline{\FF}}

\DeclareMathOperator{\Gal}{Gal}
\newcommand{\PGL}{\mathrm{PGL}}
\newcommand{\PSL}{\mathrm{PSL}}
\newcommand{\GL}{\mathrm{GL}}

\newcommand{\PSp}{\mathrm{PSp}}
\DeclareMathOperator{\Image}{im}
\newcommand{\calO}{\mathcal{O}}
\DeclareMathOperator{\Frob}{Frob}
\DeclareMathOperator{\Tr}{Tr}
\newcommand{\rhobar}{\overline{\rho}}
\newcommand{\proj}{\mathrm{proj}}
\newcommand{\mat}[4]{
 \left(  \begin{smallmatrix} #1 & #2 \\ #3 & #4 \end{smallmatrix} \right)}
\newcommand{\Ind}{{\rm Ind}}

\begin{document}

\selectlanguage{british}

{\Large Applying modular Galois representations to the Inverse Galois Problem}\\[.2cm]
{\large Gabor Wiese\footnote{Universit\'e du Luxembourg,
Facult\'e des Sciences, de la Technologie et de la Communication,
6, rue Richard Coudenhove-Kalergi, L-1359 Luxembourg, Luxembourg, gabor.wiese@uni.lu}},
5 February 2014\\[.3cm]

For many finite groups the Inverse Galois Problem (IGP) can be approached through modular/automorphic
Galois representations. This report is about the {\bf ideas} and the {\bf methods} that my coauthors and
I have used so far, and their {\bf limitations} (in my experience).

In this report I will mostly stick to the case of $2$-dimensional Galois representations because it is technically much simpler and already exhibits essential features; occasionally I'll mention $n$-dimensional
symplectic representations; details on that case can be found in Sara Arias-de-Reyna's report
on our joint work with Dieulefait and Shin.

\subsection*{Basics of the approach}

{\bf The link between the IGP and Galois representations.}
Let $K/\QQ$ be a finite Galois extension such that $G := \Gal(K/\QQ) \subset \GL_n(\Fbar_\ell)$ is a subgroup.
Then $G_\QQ := \Gal(\Qbar/\QQ) \twoheadrightarrow \Gal(K/\QQ) \hookrightarrow \GL_n(\Fbar_\ell)$ is an $n$-dimensional
continuous Galois representation with image~$G$.
Conversely, given a Galois representation $\rho: G_\QQ \to \GL_n(\Fbar_\ell)$ (all our Galois representations
are assumed continuous), then $\Image(\rho) \subset \GL_n(\Fbar_\ell)$ is the Galois group
of the Galois extension $\overline{\QQ}^{\ker(\rho)}/\QQ$.
\smallskip

{\bf Source of Galois representations: abelian varieties.}
Let $A$ be a $\GL_2$-type abelian variety over~$\QQ$ of dimension~$d$ with multiplication
by the number field $F/\QQ$ (of degree~$d$) with integer ring $\calO_F$. Then for every
prime ideal $\lambda \lhd \calO_F$, the $\lambda$-adic Tate module of~$A$ gives rise to
$\rho_{A,\lambda}: G_\QQ \to \GL_2(\calO_{F,\lambda})$.
These representations are a special case of those presented next (due to work of Ribet
and the proof of Serre's modularity conjecture).
\smallskip

{\bf Source of Galois representations: modular/automorphic forms.}
Let $f=\sum_{n=1}^\infty a_n e^{2\pi i n z}$ be a normalised Hecke eigenform of level~$N$
and weight~$k$ without CM
(or, more generally, an automorphic representation of a certain type over~$\QQ$).
The coefficients $a_n$ are algebraic integers and $\QQ_f = \QQ(a_n \;|\; n \in \NN)$ is a number field,
the {\bf coefficient field} of~$f$. Denote by $\ZZ_f$ its ring of integers.
The eigenform $f$ gives rise to a {\bf compatible system of Galois representations},
that is, for every prime~$\lambda$ of $\QQ_f$ a Galois representation
$\rho_{f,\lambda}: G_\QQ \to \GL_2(\ZZ_{f,\lambda})$
such that $\rho_{f,\lambda}$ is unramified outside $N\ell$ (where
$(\ell) = \ZZ\cap \lambda$) and for all $p \nmid N\ell$ we have
$\Tr(\rho_{f,\lambda}(\Frob_p)) = a_p$. All representations thus obtained are odd
(determinant of complex conjugation equals~$-1$).
\smallskip

{\bf Reduction and projectivisation.}
We consider the representations
$\rhobar_{f,\lambda}: G_\QQ \xrightarrow{\rho_{f,\lambda}} \GL_2(\ZZ_{f,\lambda}) \twoheadrightarrow \GL_2(\FF_{f,\lambda})$
and 
$\rhobar_{f,\lambda}^\proj: G_\QQ \xrightarrow{\rhobar_{f,\lambda}} \GL_2(\FF_{f,\lambda}) \twoheadrightarrow \PGL_2(\FF_{f,\lambda})$, where $\FF_{f,\lambda}=\ZZ_{f,\lambda}/\lambda$.
In our research we focus on projective representations because the groups $\PSL_2(\FF_{\ell^d})$
are simple for $\ell^d\ge 4$.

{\bf Main idea: By varying $f$ and $\lambda$ (and thus~$\ell$), realise as many finite subgroups of
$\PGL_2(\Fbar_\ell)$ as possible.}
\smallskip

{\bf Trust in the approach.}
If $\ell>2$, the oddness of the representations leads to $\Qbar^{\ker(\rhobar_{f,\ell}^\proj)}$
being totally imaginary.
{\bf The approach through modular Galois representations for the groups $\PSL_2(\FF_{\ell^d})$
and $\PGL_2(\FF_{\ell^d})$ to the IGP should in principle work} for
the following reason: If $\Gal(K/\QQ)\subset \PGL_2(\Fbar_\ell)$ is a finite (irreducible) subgroup
and $K/\QQ$ is totally imaginary (which is `much more likely' than being totally real), then
Serre's modularity conjecture implies that $K$ can be obtained from some $f$ and~$\lambda$.
In more general contexts, there are generalisations of Serre's modularity conjecture (however,
unproved!) and I am inclined to believe that the approach is promising in more general contexts
than just~$\GL_2$.
\smallskip

{\bf The two directions.}
We have so far explored two directions for the
realisation of $\PSL_2(\FF_{\ell^d})$ and $\PSp_n(\FF_{\ell^d})$.
{\bf Vertical direction:} fix~$\ell$, let $d$ run (results by me for $\PSL_2$~\cite{SL2},
generalised by Khare-Larsen-Savin for $\PSp_n$~\cite{KLS1});
{\bf horizontal direction:} fix~$d$, let $\ell$ run (results by Dieulefait and me for $\PSL_2$~\cite{DiWi}
and by Arias-de-Reyna, Dieulefait, Shin and me for $\PSp_n$~\cite{partIII}).

\subsection*{Main challenges}

In approaching the IGP through modular forms for specific groups,
in my experience one is faced with two challenges:\\
\hspace*{1cm} {\bf (1) Control/predetermine the type of the image $\rhobar_{f,\lambda}^\proj(G_\QQ)$.}\\
\hspace*{1cm} {\bf (2) Control/predetermine the coefficient field $\QQ_f$.}\\
Problem~(2) appears harder to me.
\smallskip

{\bf Controlling the type of the images.}
By a classical theorem of Dickson, if $\rhobar_{f,\lambda}$ is irreducible, then it is either induced
from a lower dimensional representation (only possiblity: a character)
or $\rhobar_{f,\lambda}^\proj(G_\QQ) \in \{\PSL_2(\FF_{\ell^d}),\PGL_2(\FF_{\ell^d}) \}$
for some~$d$ (we call this case {\bf huge/big image}).
Under the assumption of a transvection in the image, we have generalised this result
to symplectic representations.
In our applications we want to exclude reducibility and induction.
One can expect a {\bf generic huge image result} (for $\GL_2$ this is classical work of Ribet;
for other cases e.g.\ recent work of Larsen and Chin Yin Hui in this direction~\cite{HL}).
\smallskip

{\bf Inner twists.}
If one has e.g.\ determined that $\rhobar_{f,\lambda}^\proj(G_\QQ)$ is huge, one still needs
to compute which $d\in \NN$ and which of the two cases $\PSL_2(\FF_{\ell^d})$, $\PGL_2(\FF_{\ell^d})$ occurs.
The answer is given by {\bf inner twists}. For $\GL_2$ these are well-understood (with Dieulefait
we exclude them by a good choice of~$f$); for $\PSp_n$ we proved a generalisation allowing us
to describe~$d$ by means of a number field,
but, as to now we are unable to distinguish between the two cases.

\subsection*{Coefficient field}

One knows that $\QQ_f$ is either totally real or totally imaginary (depending on the nebentype of~$f$).
Moreover, $[\QQ_f:\QQ] \le \dim S_k(N)$, where $S_k(N)$ is the space of cusp forms of level~$N$ and
weight~$k$. Furthermore, a result of Serre says that for any sequence $(N_n,k_n)_n$ such that $N_n+k_n$
tends to infinity, there is $f_n \in S_{k_n}(N_n)$ such that $[\QQ_{f_n}:\QQ]$ tends to infinity.
However, to the best of my knowledge, almost {\bf nothing is known about the arithmetic of
the coefficient fields and the Galois groups of their normal closures over~$\QQ$}.
In my experience, this is the {\bf biggest obstacle} preventing us from obtaining very strong results on
the IGP.
\smallskip

{\bf Almost complete control through Maeda's conjecture.}
A conjecture of Maeda gives us some control on the coefficient field by claiming
that for any $f \in S_k(1)$ one has $[\QQ_f:\QQ]=\dim S_k(1) =: m_k$
and that the Galois group of the normal closure of $\QQ_f$ over~$\QQ$ is $S_{m_k}$, the symmetric group.
The conjecture has been numerically tested for quite high values of~$k$, but to my knowledge a proof
is out of sight at the moment and there's no generalisation to higher dimensions either.
Assuming Maeda's conjecture I was able to prove in~\cite{Maeda} that for even~$d$
the groups $\PSL_2(\FF_{\ell^d})$ occur as Galois groups over~$\QQ$ with only $\ell$ ramifying
for all~$\ell$, except possibly a density-$0$ set. In a nutshell, for the proof I
choose a sequence $f_n$ of forms of level~$1$ such that $[\QQ_{f_n}:\QQ]$ strictly increases.
That the Galois group is the symmetric group ensures two things: firstly, every $\QQ_{f_n}$ possesses
a degree-$d$ prime; secondly, the fields $\QQ_{f_n}$ and $\QQ_{f_m}$ for $m \neq n$ are almost
disjoint (in the sense that their intersection is at most quadratic) and thus the sets of primes of degree~$d$
in the two fields are almost independent, so that their density adds up to~$1$ when $n \to \infty$.
This illustrates that {\bf some control on the coefficient field promises strong results on the IGP}.
\smallskip

{\bf A conjecture of Coleman on $\GL_2$-type abelian varieties.}
The modular form~$f$ corresponding to a $\GL_2$-type abelian variety with multiplication by~$F$
has coefficient field $\QQ_f=F$. However, I don't know of any method to construct a
$\GL_2$-type abelian variety with multiplication by a given field.
Indeed, a conjecture attributed to Coleman (see~\cite{BFGR}) predicts
that for a given dimension, only finitely many number fields occur. In other words, for weight-$2$
modular forms in all levels, there are only finitely many $\QQ_f$ of a given degree.
Under the assumption of Coleman's conjecture, it is impossible to obtain $\PSL_2(\FF_{\ell^2})$
for all~$\ell$ from $\GL_2$-type abelian surfaces because there will be a positive density set
of~$\ell$ that are split in all number fields of degree~$2$ that occur as multiplication fields.
Although I don't know if there are finitely or infinitely many quadratic fields occuring as $\QQ_f$ for
$f$ of arbitrary level and arbitrary weight, this nevertheless suggests to me that one should make use
of modular forms of {\bf arbitrary coefficient degrees} for approaching $\PSL_2(\FF_{\ell^d})$ for fixed~$d$
(as we did when we assumed Maeda's conjecture).
\smallskip

{\bf Numerical data.}
Some very simple computer calculations for $p=2$ during my PhD have very quickly revealed that all
$\PSL_2(\FF_{2^d})$ with $1\le d\le 77$ occur over~$\QQ$. With Marcel Mohyla we
plotted $\FF_{f,\lambda}$ for small fixed weight and $f$ having prime levels~\cite{MaGa}. The
computations suggest that the maximum and the average degrees (for $f$ in $S_k(N)$ for $N$ prime)
of $\FF_{f,\lambda}$ are roughly proportional to the dimension of~$S_k(N)$.

\subsection*{The local `bad primes' approach to the main challenges}

We need to gain some control on the coefficient fields and in the absence of a generic huge image
result, we also need to force huge image of the Galois representation. In all our work
(like in that of Khare-Larsen-Savin~\cite{KLS1}),
we approach this by choosing suitable inertial types, or in the
language of abelian varieties, by choosing certain types of bad reduction. The basic idea appeared
in the work of Khare-Wintenberger on Serre's modularity conjecture. More precisely, one chooses
inertial types at some primes~$q$ guaranteeing that $\rhobar_{f,\lambda}(I_q)$ contains certain elements
($I_q$ denotes the inertia group at~$q$) .
For instance, if an element that is conjugate to $\mat 1101$ is contained,
the representation cannot be induced. In the $n$-dimensional symplectic case, we use this to
obtain a transvection in the image, allowing us to apply our classification (see above). We also
employ Khare-Larsen-Savin's generalisation of Khare-Wintenberger's good-dihedral primes. More precisely,
for $\GL_2$ we impose $\rhobar_{f,\lambda}|G_{\QQ_q} = \Ind_{\QQ_{q^2}}^{\QQ_q}(\alpha)$ where $\alpha$
is a character of $\QQ_{q^2}^\times$ of prime order~$t$ not descending to $\QQ_q^\times$.
This has two uses: (1) As the representation is {\bf irreducible} locally at~$q$, so it is globally.
(2) $\QQ_f$ contains $\zeta_t+\zeta_t^{-1}$ (this follows from an explicit description of the induction).
This {\bf cyclotomic field in the coefficient field} can be exploited in two ways.
(2a) By making $t$ big, $[\FF_{f,\lambda}:\FF_\ell]$ becomes big. {\bf This leads to the results
in the vertical direction.}
(2b) Given~$d$, by choosing $t$ suitably, $\QQ(\zeta_t+\zeta_t^{-1})$ contains prime ideals of degree~$d$,
thus {\bf $\QQ_f$ contains prime ideals of degree~$d$, which makes the results in the horizontal
direction work}.
In the absence of any knowledge on the Galois closure of $\QQ_f$ over~$\QQ$ in general, I do not know of
any other way to guarantee that degree-$d$ primes exist at all (we need them to realise
$\PSL_2(\FF_{\ell^d})$).

My feeling is that the cyclotomic field $\QQ(\zeta_t+\zeta_t^{-1})$ only makes up a very small
part of the coefficient field, i.e.\ that $[\QQ_f:\QQ]$ will be much bigger than
$[\QQ(\zeta_t+\zeta_t^{-1}):\QQ]$. Thus, in our results in the horizontal direction, for given $d$ and~$f$,
we only obtain very small densities.
Moreover, I cannot prove that by varying~$f$ for fixed~$d$,
the sets of primes of residue degree~$d$ are not contained in each other. Any information, for instance,
on the ramification of $\QQ_f$ changing with~$f$ or on the Galois group would probably enable us to
obtain a big density by taking the union of the sets of degree-$d$ primes for many~$f$.

\subsection*{Constructing the relevant modular/automorphic forms}

For finishing the approach, one must finally construct or show the existence of modular/automorphic
forms having the required inertial types. For modular forms one can do this in quite a down-to-earth
way by using level raising. This approach was taken in the work by Dieulefait and me.
In the symplectic case, we exploit work of Shin, as well as level-lowering results of
Barnet-Lamb, Gee, Geraghty and Taylor~\cite{BLGGT}.
Khare-Larsen-Savin~\cite{KLS1} use other automorphic techniques.

\subsection*{Conclusion}
The presented approach to the IGP for many families of finite groups
through automorphic representations seems in principle promising.
In my opinion, the main obstacle is a poor understanding of the coefficient fields.

The approach has the advantage that it allows {\bf full control on the ramification}.
A disadvantage is that one does not obtain a regular realisation.
\smallskip

{\bf Acknowledgements.} I thank Sara Arias-de-Reyna for valuable comments on a first draft
of this report.

\bibliography{Bibliog}
\bibliographystyle{amsalpha}

\end{document}